\documentclass[11pt,draft]{amsart}

\usepackage{amssymb}
\usepackage{url}
\usepackage{amscd}
\usepackage{amsthm}
\usepackage{amsmath}

\theoremstyle{plain}
\newtheorem{theorem}{Theorem}[section]

\theoremstyle{remark}

\begin{document}

\date{}

\title[Nikodym sets]
{On the size of Nikodym sets in finite fields}
\author{Liangpan Li}
\date{}
\address{Department of Mathematics,
Shanghai Jiaotong University, Shanghai 200240, People's Republic of
China} \email{liliangpan@yahoo.com.cn}

\begin{abstract}
Let $\mathbb{F}_q$ denote a finite field of $q$ elements. Define a
set $B\subset\mathbb{F}_q^n$ to be Nikodym if for each $x\in B^{c}$,
there exists a line $L$ such that $L\cap B^c=\{x\}.$ The main
purpose of this note is to show that the size of every Nikodym set
is at least $C_n\cdot q^n$, where $C_n$ depends only on $n$.
\end{abstract}

\maketitle

\section{Introduction}
The finite field Kakeya problem, posed by Wolff in his influential
survey \cite{wolff}, asks for the smallest subset of
$\mathbb{F}_q^n$ that contains a line in each direction, where
$\mathbb{F}_q$ denotes a finite field of $q$ elements. A subset
containing a line in each direction is called a Kakeya set. In
analogy with the Euclidean Kakeya problem, Wolff conjectured that
 $\sharp K\geq C_{n}q^n$
holds for any Kakeya set $K\subset\mathbb{F}_q^n$, where $C_n$
depends only on the dimension $n$. For $n=2$ Wolff immediately
proved the bound $\sharp K\geq q(q+1)/2$, and it is best possible
when $q$ is even. To the author's knowledge, Blokhuis and Mazzocca
\cite{blokhuis} studied the finite field Kakeya problem in two
dimensions and proved the sharp bound $\sharp K\geq
{q(q+1)}/{2}+({q-1})/{2}$ when $q$ is odd, as conjectured by Faber
in \cite{faber}. The higher dimensional finite field Kakeya problem
has been extensively investigated
 in
\cite{bourgainkatztao,mockenhaupttao,rogers,taofourdimension,wolff}
such as proving the bound $\sharp K\geq C_nq^{(n+2)/2}$ or $\sharp
K\geq C_nq^{(4n+3)/7}$. Recently, using the polynomial method in
algebraic extremal combinatorics, Dvir \cite{Dvir} completely
confirmed this conjecture by proving
\[\sharp K\geq{n+q-1\choose n}.\]




On the other hand, Nikodym \cite{nikodym} proved that there exists a
null set in the unit square such that every point of the complement
is ``linearly accessible through the set", which means it lies on a
line that is otherwise included in the set. Falconer \cite{falconer}
extended Nikodym's result to higher dimensions proving there exists
a set $N\subset\mathbb{R}^n$ of zero Lebesgue measure such that for
each $x\in N^c$, there is a hyperplane $P$ satisfying $P\cap
N^c=\{x\}$. In the Euclidean spaces Nikodym sets are closely related
to Kakeya sets through Carbery's transformation \cite{Carbery,tao}.

Motivated by the above works,   we shall define  a set $B$ in
$\mathbb{F}_{q}^n$ to be \textsf{Nikodym} if for each $x\in B^{c}$
there exists a line $L$ such that $L\cap B^c=\{x\}$. The main
purpose of this note is to prove the lower bound
\[\sharp B\geq{n+q-2\choose n}.\]
Slightly different with the two dimensional finite field Kakeya
problem, this bound is  not best possible in two dimensions.





\section{General dimensions}

\begin{theorem}
 Any Nikodym
set $B\subset\mathbb{F}_{q}^n$  satisfies
\[|B|\geq{n+q-2\choose n},\]
where $\mathbb{F}_{q}$ denotes a finite field of $q$ elements.
\end{theorem}



\begin{proof}
We argue by contradiction and suppose
\[|B|<{n+q-2\choose n}.\]
A basic result in combinatorics \cite{Brualdi} says that the number
of monomials in $\mathbb{F}[x_1,\ldots,x_n]$ of degree at most $d$
is
\[{{n+d \choose n}},\]
hence there exists a nonzero polynomial
$g\in\mathbb{F}[x_1,\ldots,x_n]$ of degree at most $q-2$ such that
\[g(y)=0 \ \ (\forall y\in B).\] For each $x\in B^{c}$, there exists a
line $L$ such that
\[L\cap B^c=\{x\}.\]
The restriction of $g$ to this line is a univariate polynomial of
degree at most $q-2$, and since it has at least $q-1$ zeros, it must
be zero on the entire line $L$. Considering $x$ belongs to this
line, it follows that
\[g(x)=0.\]
This would mean $g$ is the zero polynomial, a contradiction.

\end{proof}

\section{Two dimensions }


\begin{theorem}
 Any Nikodym
set $B\subset\mathbb{F}_{q}^2$  satisfies
\[\sharp B\geq\frac{2q^2}{3}+O(q) \ \ (q\rightarrow\infty),\]
where $\mathbb{F}_{q}$ denotes a finite field of $q$ elements.
\end{theorem}

\begin{proof}
Write $s=\lfloor\frac{q}{3}\rfloor$. First,  assume that
\[\sharp
B^c\leq s(q-1)+2q,\]
then \begin{equation}
\sharp B\geq
q^2-s(q-1)-2q\geq
q^2-\frac{q}{3}(q-1)-2q=\frac{2q^2}{3}-\frac{5q}{3}.\label{estimate
1}
\end{equation}
Else suppose that
\[\sharp B^c\geq s(q-1)+2q.\]
Since $B$ is a Nikodym set, for each $x\in B^c$ there exists a line
$L_x$ such that
\[L_x\cap B^c=\{x\}.\]
Obviously, all of these lines are distinct from each other.
 Noting that there are in total $q+1$ directions in
$\mathbb{F}_{q}^2$, we partition $\{L_x\}_{x\in B^c}$ into classes
$\{G_{i}\}_{i=0}^{q}$ according to their directions. Without loss of
generality we may assume that
\[\sharp G_0\geq\sharp G_1\geq\sharp
G_2 \geq\cdots\geq\sharp G_q.\]
 Thus
\[q+q+\sharp G_2\cdot(q-1)\geq\sum_{i=0}^{q}\sharp G_i=\sharp B^c\geq s(q-1)+2q,\]
from which yields
\[\sharp G_2\geq s.\]
Choose $s$ parallel lines $\{X_l\}_{l=1}^{s}$ from $G_0$, $s$
parallel lines $\{Y_m\}_{m=1}^{s}$ from $G_1$ and $s$ parallel lines
$\{Z_n\}_{n=1}^{s}$ from $G_2$, then it follows that
\begin{align}
\sharp B&\geq
\sum_{l=1}^{s}(\sharp X_l-1)+\sum_{m=1}^{s}(\sharp
Y_m-1-s)+\sum_{n=1}^{s}(\sharp Z_n-1-2s)\nonumber\\
&=s(q-1)+s(q-1-s)+s(q-1-2s)=3s(q-1-s)\nonumber\\
&\geq3\frac{q-2}{3}(q-1-\frac{q}{3})=\frac{2q^2}{3}-\frac{7q}{3}+2.\label{estimate
2}
\end{align}
Combining (\ref{estimate 1}) and (\ref{estimate 2}) yields the
desired result.

\end{proof}

\textsc{Question}: How small can the Nikodym sets really be in two
dimensions?


\section{Acknowledgements}
The author thanks  Yaokun Wu for clarifying the proof. He also
thanks Aart Blokhuis and  Qing Xiang for kindly pointing out the
recent progresses on the finite field Kakeya problem to the author.

\end{document}